\documentclass{article}

\usepackage{arxiv}

\usepackage[utf8]{inputenc} % allow utf-8 input
\usepackage[T1]{fontenc}    % use 8-bit T1 fonts
\usepackage{hyperref}       % hyperlinks
\usepackage{url}            % simple URL typesetting
\usepackage{booktabs}       % professional-quality tables
\usepackage{amsfonts}       % blackboard math symbols
\usepackage{nicefrac}       % compact symbols for 1/2, etc.
\usepackage{microtype}      % microtypography
\usepackage{lipsum}
\usepackage{amssymb}        %symbols
\usepackage{amsmath}        %equations

\title{Infinite Sequences, Series Convergence and the Discrete Time Fourier Transform over Finite Fields}

\author{
  R. M.~Campello~de~Souza \\
  Federal University of Pernambuco\\
  50711-970, Recife, PE Brazil.\\
  \texttt{ricardo@ufpe.br} \\
   \And
 M. M.~Campello~de~Souza \\
  Federal University of Pernambuco\\
  50711-970, Recife, PE Brazil.\\
  \texttt{marciam@ufpe.br} \\
  \And
  H. M.~de~Oliveira\thanks{\texttt{https://arxiv.org/a/deoliveira\_h\_1.html}~---~\texttt{https://orcid.org/0000-0002-6843-0635.}}\\
  Federal University of Pernambuco\\
  50711-970, Recife, PE Brazil.\\
  \texttt{hmo@ufpe.br} \\
  \And
  M. M.~Vasconcelos\\
   Federal University of Pernambuco\\
  50711-970, Recife, PE Brazil.\\
  \texttt{mmv@ee.ufpe.br} \\
}

\begin{document}
\maketitle

\begin{abstract}
Digital Transforms have important applications on subjects such as channel coding, cryptography and digital signal processing. In this paper, two Fourier Transforms are considered, the discrete time Fourier transform (DTFT) and the finite field Fourier transform (FFFT). A finite field version of the DTFT is introduced and the FFFT is redefined with a complex kernel, which makes it a more appropriate finite field version of the Discrete Fourier Transform. These transforms can handle FIR and IIR filters defined over finite algebraic structures.
\end{abstract}

% keywords can be removed

\section{Introduction}
Discrete Transforms play a very important role in engineering. Well known examples are the Discrete Fourier Transform (DFT) and the Z Transform [\cite{Oppenheim}]. A DFT over finite fields GF(q) (FFFT) was also defined \cite{Pollard} and applied to the computation of discrete convolutions through modular arithmetic. Recently, the Hartley Transform over $GF(q)$ (FFHT) was introduced \cite{Campello}, which has interesting applications in digital multiplexing and spread spectrum [4, 5]. The FFFT and the FFHT, instead of what happens with others discrete transforms, are examples of transforms of true digital nature.
In this paper, the Discrete Time Fourier Transform over a finite field (FFDTFT) is introduced. So far, the investigations reported in the literature about transforms defined over finite fields, only consider finite sequences of elements from $GF(q)$, a constraint related to the problem of dealing with series convergence over a finite algebraic structure \cite{Cooklev}. The proposed FFDTFT deals with finite and infinite sequences over finite fields. The paper also introduces a generalization of the DFT over a finite field, which results from the modification of its kernel. An FFFT with complex kernel is the appropriate finite field version for the DFT. In the next section, some mathematical preliminaries are presented. In particular, the Gaussian integers over $GF(p)$ (denominated Galoisian integers) are defined and some special group families are built, in order to construct a polar representation for the Galoisian integers. In Section 3, infinite sequences over a Galois field are discussed, where issues on series convergence are investigated. In Section 4, the Discrete Time Fourier Transform over a finite field is introduced. In Section 5, a new definition for the FFFT is proposed. The paper closes with a few concluding remarks.

\section{Mathematical Preliminaries}
\label{sec:Preliminaries}

\subsection{Complex Numbers over Finite Fields}
The set $Ga(p)$ of Gaussian integers over $GF(p)$ defined below plays an important role in the ideas introduced in this paper. Hereafter the symbol $\triangleq$ denotes \textit{equal by definition}; $\mathbb{Q}$, $\mathbb{R}$ and $\mathbb{C}$ denote the rational, real and complex sets, respectively. $\delta$ is the Kronecker symbol.

\textit{\textbf{Definition 1}. $Ga(p) \triangleq \{a + jb; ~a, b \in GF(p)\}$, $p$ being an odd prime for which $j^2=-1$ is a quadratic non-residue in $GF(p)$ (i.e., $p \equiv 3 \pmod{4}$), is the set of Galoisian integers over $GF(p)$}. \\

Let $\otimes$ denote Cartesian product. It can be shown that the set $Ga(p)$ equipped with the operations $\oplus$ and $*$ defined below, is a field \cite{Blahut}.\\

\textit{\textbf{Proposition 1}. Let}
\begin{align*}
&\oplus :Ga(p) \otimes Ga(p) \to  Ga(p)\\
&(a_1+jb_1,a_2+jb_2) \to 
(a_1 +jb_1) \oplus (a_2 +jb_2) = (a_1 +a_2)+j(b_1+b_2)
\end{align*}
\\\textit{and}\\
\begin{align*}
&*: Ga(p) \otimes Ga(p) \to  Ga(p)\\ 
\\
&(a_1+jb_1,a_2+jb_2) \to  (a_1+jb_1)*(a_2+jb_2) =
(a_1a_2-b_1b_2) + j(a_1b_2+a_2b_1).
\end{align*}

\textit{The structure $GL(p) \triangleq <Ga(p), \oplus ,* >$ is a field. In fact, $GL(p)$ is isomorphic to $GF(p^2)$.}\\

By analogy with the real and complex numbers, the elements of $GF(p)$ and of $GL(p)$ are said to be real and complex, respectively.\\

\textit{\textbf{Proposition 2.}} \textit{The elements $\zeta=(a+jb) \in GL(p)$ satisfies $\zeta ^{p+1} \equiv |\zeta|^2 \equiv a^2 + b^2 \pmod{p}$}.\\

\textit{Proof}:\\
Once $GL(p)$ is isomorphic to $GF(p^2)$, a field of characteristic $p$. Since $p = 4k + 3$, $jp = -j$, so that $\zeta^{p} \equiv a-jb \pmod{p} = \zeta^* \pmod{p}$. Therefore, $\zeta^{p+1} \equiv \zeta \zeta^* = |\zeta|^2 \equiv a^2+b^2 \pmod{p}$.

\subsection{Polar Form for Galoisian Integers in a Finite Field}

In the definition of $GL(p)$, the elements were written in Cartesian form $\zeta = a + jb$. In what follows, a different representation for the elements of the multiplicative group of $GL(p)$ is proposed, which allows to write them in the form $r\epsilon ^\theta$. By analogy with the continuum, such a form is said to be polar.\\

\textit{\textbf{Proposition 3}. Let $G_r$ and $G_{\theta}$ be subgroups, of the multiplicative group $G$ of the nonzero elements of $GL(p)$, of orders $N_r = (p-1)/2$ and $N_{\theta} = 2(p+1)$, respectively. Then all elements of $GL(p)$ can be written in the form $\zeta=ab$, where $a \in G_r$ and $b \in G_{\theta}$} \cite{Lidl}.\\

Considering that any element of a cyclic group can be written as an integer power of a group generator, it is possible to set $r=a$ and $\epsilon ^\theta =b$, where $\epsilon$ is a generator of $G_{\theta}$. Thus, the polar representation assumes the desired form, $\zeta = r \epsilon ^\theta$.\\

At this point, it seems clear that $r$ is going to play the role of the modulus of $\zeta$. Therefore, before further exploring the polar notation, it is necessary to formally define the concept of modulus of an element in a finite field. Considering the nonzero elements of $GF(p)$, it is a well-known fact that half of them are quadratic residues of $p$ \cite{Burton}. The other half, those that do not have a square root, are the quadratic non-residues. Likewise, in the field $\mathbb{R}$ of real numbers, the elements are divided into positive and negative numbers, those that have and those that do not have a square root. The standard modulus operation in $\mathbb{R}$ always gives a positive result. By analogy, the modulus operation in $GF(p)$ is going to be defined, such that it always results in a quadratic residue of $p$.\\

\textit{\textbf{Definition 2}. The modulus of an element $a \in  GF(p)$, where $p = 4k + 3$, is given by}
\begin{equation*}
\left | a \right | \triangleq \left\{\begin{matrix}
a & if ~~a^{(p-1)/2} \equiv 1 \pmod{p}\\ 
-a & if ~~a^{(p-1)/2} \equiv -1 \pmod{p}.
\end{matrix}\right.
\end{equation*}

\textit{\textbf{Proposition 4}. The modulus of an element of $GF(p)$ is a quadratic residue of $p$}.\\

\textit{Proof}: Since $p=4k+3$, it implies that $(p-1)/2-2k+1$, such that $(-1)^{(p-1)/2} \equiv -1 \pmod{p}$. By Euler’s criterion \cite{Burton}, if $a(p-1)/2 \equiv 1 \pmod{p}$, then $a$ is a quadratic residue of $p$; if $a^{(p-1)/2} \equiv -1 \pmod{p}$, then $a$ is a quadratic non-residue of $p$. Therefore, $(-a)^{(p-1)/2} \equiv (-1)(-1) \equiv 1 \pmod{p}$ and it follows that $|a|$ is a quadratic residue of $p$.\\

\textit{\textbf{Definition 3}. The modulus of an element $a + jb \in GL(p)$, where $p=4k+3$, is given by}

\begin{align*}
\left | a+jb \right | \triangleq  \left | \sqrt{\left | a^2+b^2 \right |} \right |.
\end{align*}

The inner modulus sign in the above expression is necessary in order to allow the computation of the square root of the quadratic norm $a^2 + b^2$, and the outer one guarantees that such an operation results in one value only. In the continuum, such an expression reduces to the usual norm of a complex number, since both, $a^2 + b^2$ and the square root operation, produce only positive numbers.\\

\textit{\textbf{Proposition 5}. If $\zeta = a+jb = r\epsilon ^\theta$, where $r \in G_r$ and $\epsilon ^\theta \in G_{\theta}$, then $r=|\zeta|$}.\\

\textit{Proof}: Every element of $G_r$ has an order that divides $(p-1)/2$. Thus, if $r \in  G_r$, then $r(p-1)/2 \equiv 1 \pmod{p}$, and $|r| = r$. Every element $\gamma$ of $G_{\theta}$ has order that divides $2(p + 1)$ and are those $c + jd$ satisfying $c^2+d^2 \equiv \pm 1 \pmod{p}$, since that, from Proposition 2, $\gamma^{2(p+1)} \equiv (c^2+d^2)^2 \equiv 1 \pmod{p}$. Besides that, as shown in next section, the elements of the group $G_{\theta}$ are those $a+jb$ such that $a^2 +b^2 \equiv \pm 1 \pmod{p}$. Therefore, according to Definition 2, such elements have modulus equal to one, which means that $|\zeta|=|r \epsilon ^\theta|=|r||\epsilon^\theta|=r.1=r$.\\

From the above it can be observed that the polar representation being introduced is consistent with the usual polar form defined over the complex field $\mathbb{C}$. The modulus belongs to $GF(p)$ (the modulus is a real number) and is a quadratic residue of $p$ (a positive number), and the exponential component $\epsilon^\theta$ has modulus one and belongs to $GL(p)$ ($e^{j\theta}$ also has modulus one and belongs to the complex field).

\section{Infinite Series over Finite Fields}

Given a sequence of integers $ \left \{ x[n] \right \} _{-\infty }^{+\infty }$, it is possible to generate a sequence over $GF(p)$ simply by considering $  \left \{ x[n] \pmod{p} \right \} _{-\infty }^{+\infty }$. In general, $x[n]$ may even be a sequence of rational elements and any element $r/s \in \mathbb{Q}$ can be mapped over $GF(p)$, $[r \pmod{p}].[s \pmod{p}]^{-1}$. The focus here concerns finite sequences, and periodic infinite sequences over a finite field.\\

\textit{\textbf{Definition 4}. An infinite sequence of elements over $GF(p)$ is periodic with period $P$, if it satisfies $x[n] = x[n \pmod{p}]$}.\\

\textit{\textbf{Definition 5}. An infinite sequence, of elements over $GF(p)$, that is zero valued for $n < N_1 < \infty$ ($-\infty < N_2 < n$) and satisfies the condition on Definition 4, is denominated left-sided periodic (right-sided periodic, respectively)}.\\
\vskip 1.5 cm
For instance, consider the following right-sided sequences:

\textbf{Example 1}. Let $ \left \{ x[n] \right \} _{-\infty }^{+\infty }$, $x[n] \in GF(p)$, e.g.
\begin{itemize}
\item  $ \left \{ 3^n \right \} _{0}^{+\infty }$ (in $GF(7)$)= $\{1~3~2~6~4~5~1~3~2~6~4~5 \cdots\}$ (P = 6).
\item $ \left \{ 1^n \right \} _{0}^{+\infty }$ (in $GF(5)$)= $\{1~1~1~1~1~1~1~1~1~1~1~1 \cdots\}$ (P = 1). 
\item $ \left \{ x[n] \right \} _{-\infty}^{+\infty }$ (in $GF(3)$)= $\{\cdots 0~0~0~1~1~1~2~2~2~0~0~0 \cdots\}$ (P = 9).
\end{itemize}

Is it possible to define a Discrete Time Fourier Transform for these sequences, by $X(\epsilon^\theta) \triangleq \sum_{n=-\infty }^{+\infty }x[n]\epsilon ^{-n\theta}$
where $\epsilon \in GL(p)$? Is this series convergent?\\

It is a well-known fact that the infinite series
$\sum_{n=1 }^{+\infty }(-1)^{n+1}= 1-1+1-1+1-\cdots$ diverges in the classic sense. However, Euler among others noticed that the arithmetic mean of the partial sums converges to 1/2. The partial sums of this series are $S_1 =1,~S_2 =0,~S_3 =1,~S_4 =0,\cdots$ and the arithmetic mean $\sigma_n \triangleq (1/n) \sum_{k=1}^n S_k$ forms the sequence ($\sigma_n$) that converges to 1/2.
When a series converges in the sense that the arithmetic mean of the partial sums converges, it is said to be Ces\`aro-summable (Ernesto Ces\`aro (1859- 1906)) \cite{Bartle}. Every convergent series in the usual sense is Ces\`aro-summable and the series sum is equal to the limit of the sequence of the partial sums arithmetic mean. That shows the Ces\`aro summability, is now introduced. Given $ \left \{ x[n] \right \} _{1}^{+\infty }$, the partial sums $S[n]$ are defined according to: 
$S[n] \triangleq \sum_{k=1 }^{n}x[k]$.\\

\textit{\textbf{Definition 6}. The Cesa\`ro sum over a finite field is defined by}
\begin{equation*}
\sigma_n \triangleq \frac{1}{n} \sum_{k=1}^n S[k],
\end{equation*}
\textit{where $S[k] \in GF(p)$ are interpreted as integers}.\\

If $ \left \{ x[n] \right \} _{1}^{+\infty }$ is a periodic sequence over $GF(p)$, so is $ \left \{ S[n] \right \} _{1}^{+\infty }$. Let $P$ denote the period of the latter sequence. Therefore
\begin{equation*}
\lim_{n \to \infty }\frac{1}{n}\sum_{k=1}^{n} S[k]=\lim_{n \to \infty } \left ( \sum_{k=1}^{\left \lfloor n/P \right \rfloor P} S[k]  + \sum_{k=1 \pmod{p})}^{n} S[k] \right )
\end{equation*}

The second term, which exists only if $P$ do not divide $n$, vanishes and
\begin{equation*}
\lim_{n \to \infty } \sigma_n=\lim_{n \to \infty }\frac {\lfloor n/P \rfloor}{n} \sum_{k=1}^P S[k].
\end{equation*}
But
\begin{equation*}
\lim_{n \to \infty }\frac {\lfloor n/P \rfloor}{n}= \frac{1}{P},
\end{equation*}
so that 
\begin{equation*}
\lim_{n \to \infty }\frac {1}{n} \sum_{k=1}^n S[k]=\frac {1}{P} \sum_{k=1}^P S[k],
\end{equation*}
and therefore
\begin{equation*}
\lim_{n \to \infty } \sigma_n \pmod{p} \equiv \lim_{n \to \infty }\frac {1}{P \pmod{p}} \sum_{k=1}^P S[k] \pmod{p}.
\end{equation*}
The heart of the matter is taking first the limit $n \to \infty$, and then evaluate the result after reducing it modulo $p$.\\

\textit{\textbf{Definition 7}. A series over a finite field is said to be Ces\`aro convergent to $\sigma$ if and only if}

\begin{equation*}
\sigma \triangleq \lim_{n \to \infty } \sigma_n \pmod{p} \in GF(p).
\end{equation*}

\textit{\textbf{Corollary 1}. Every periodic series over a finite field with a nonzero period, that is, $P \neq 0 \pmod{P}$, is Ces\`aro convergent}.\\

\textbf{Example 2}. Considering the sequence $ \left \{ 3^n \right \} _{0}^{+\infty }$ (in $GF(7)$) from Example 1,\\
${S[k]} = \{1~4~6~5~2~0~1~4~6~5~2~0~1~4~\cdots\}$, $P \equiv 6 \pmod{7}$. Therefore, the series converges, in the Ces\`aro sense, to
\begin{equation*}
\sigma = \frac{1}{6} (1+4+6+5+2+0) \equiv 3 \pmod{7}.
\end{equation*}

\section{The Discrete Time Fourier Transform in a Finite Field}

\subsection{Basic Sequences}

The transforms considered in this paper deal with sequences $x[n]$, defined over the finite field $GF(p)$, which are obtained from the basic sequences  $\delta[n], u[n] ~\text{and } Aa^n$.

i. The finite field impulse over $GF(p)$ (Galois impulse), denoted by $\delta[n]$, is the sequence $x[n]$ defined by\\

\begin{equation*}
x[n]=\delta[n] \triangleq \left\{\begin{matrix}
1 & if~n \equiv 0 \pmod{2(p+1))}\\ 
0, & otherwise.
\end{matrix}\right.
\end{equation*}
By analogy with sequences defined over the infinite field $\mathbb{R}$, any sequence $x[n]$ defined over a finite field can be expressed as a sum of scaled and time-shifted Galois impulses.\\

ii. The finite field unit step is given by

\begin{equation*}
x[n]=u[n] \triangleq  \left\{\begin{matrix}
1 & if~n \geq 0\\ 
0 & otherwise.
\end{matrix}\right.
\end{equation*}

iii. The exponential sequence is $x[n] =A(a)^n,~A \text{ and } a \in GF(p)$. This sequence is periodic with period $P$, which is the multiplicative order of $a\pmod{p}$.\\

\textit{\textbf{Definition 8}. The Discrete Time Fourier Transform (DTFT) of a sequence $x[n]$ over $GF(p)$ is the function $X(\epsilon^\theta)$, defined in $GL(p)$, given by}
\begin{equation*}
X(\epsilon^\theta) \triangleq \sum_{n=-\infty}^{+\infty}x[n]\epsilon ^{-n\theta},
\end{equation*}
where $\epsilon \in G_{\theta}$ has multiplicative order $2(p + 1)$.\\

In the infinite series defined above, the convergence is considered in the sense of Definition 7.\\

\textbf{Example 3}. The right exponential sequence over $GF(p)$. Let $x[n] = a^n u[n], ~a \in GF(p)$. In case, since $x[n]$ is nonzero only for $n \geq 0$, then 

\begin{equation*}
X(\epsilon^\theta) = \sum_{n=-\infty }^{+\infty} a^n u[n] \epsilon ^{-n\theta}= \sum_{n=0}^{+\infty} (a \epsilon)^n.
\end{equation*}
Computing the partial sums:
\begin{align*}
S_1&=1\\
S_2&=1+a \epsilon^{-\theta}\\
S_3&=1+a \epsilon^{-\theta}+(a \epsilon^{-\theta})^2\\
&~~\vdots ~~~~~\vdots  ~~~~~~~\ddots\\
S_{N-2}&=1+a \epsilon^{-\theta}+(a \epsilon^{-\theta})^2+\cdots+(a \epsilon^{-\theta})^{N-2}\\
S_{N-1}&=1+a \epsilon^{-\theta}+(a \epsilon^{-\theta})^2+\cdots++(a \epsilon^{-\theta})^{N-2}+(a \epsilon^{-\theta})^{N-1}\\
\end{align*}

Denoting by $N$ the multiplicative order of $(a \epsilon^{-\theta})$, the sequence $S[k]$ has a period equal to $N$. Therefore, it is possible to write
\begin{equation*}
\sigma_N=\frac{1}{N} \sum_{i=0}^{N-1} (N-i)(a \epsilon^{-\theta})^i,
\end{equation*}
or
\begin{equation*}
\sigma_N=\frac{(a \epsilon^{-\theta})^N-1}{a \epsilon^{-\theta}-1}-\frac{1}{N}\sum_{i=0}^{N-1} i(a \epsilon^{-\theta})^i.
\end{equation*}
Since $(a \epsilon^{-\theta})$ has multiplicative order $N$,
\begin{equation*}
\sigma_N=-\frac{\epsilon^{\theta}}{N}\sum_{i=0}^{N-1} i(a \epsilon^{-\theta})^{-i-1}.
\end{equation*}
which is the same as
\begin{equation*}
\sigma_N=-\frac{\epsilon^{\theta}}{N}\sum_{i=0}^{N-1} \frac{\mathrm{d} }{\mathrm{d} \epsilon^{\theta}}  (a^i (\epsilon^{-\theta})^{-i}),
\end{equation*}
so that
\begin{equation*}
\sigma_N=-\frac{\epsilon^{\theta}}{N}  \frac{\mathrm{d} }{\mathrm{d} \epsilon^{\theta}} \sum_{i=0}^{N-1} (a \epsilon^{-\theta})^i,
\end{equation*}
A simple manipulation shows that the last expression is equal to
\begin{equation*}
\sigma_N=\frac{1}{1-a \epsilon^{-\theta}}.
\end{equation*}
which is the FFDTFT $X(\epsilon^\theta)$.

\subsection{The Inverse FFDTFT}

\textit{\textbf{Lemma 1}. If $\epsilon \in G_\theta$ has multiplicative order $2(p + 1)$, then}

\begin{equation*}
\sum_{\theta=0}^{2(p+1)-1} \epsilon^{\theta k}=\left\{\begin{matrix}
2(p+1), & if ~k \equiv 0 \pmod{2p+1}\\ 
0, & otherwise.
\end{matrix}\right.
\end{equation*}

\textit{Proof}. For $k \equiv 0 \pmod{2p+1}$, the sum is clearly equal to $2(p + 1)$. Otherwise,
\begin{equation*}
\sum_{\theta=0}^{2(p+1)-1} \epsilon^{\theta k}=\frac{1-\epsilon^{k2(p+1)}}{1-\epsilon^k}
\end{equation*}
\textit{and the result follows}.\\

\textit{\textbf{Theorem 1}. (the inversion formula) The inverse finite field discrete time Fourier transform is given by}
\begin{equation*}
x[n]=\frac{1}{2(p+1)}\sum_{\theta =0}^{2(p+1)-1}X(\epsilon^\theta) \epsilon^{\theta n}.
\end{equation*}
\textit{Proof}. By definition $X(\epsilon^{\theta})=\sum_{k=-\infty }^{+\infty} x[k]\epsilon^{k\theta}$. Multiplying both sides by $\epsilon^{n\theta}$ and summing over $\theta$, we have
\begin{equation*}
\sum_{\theta =0}^{2(p+1)-1} X(\epsilon^\theta)\epsilon^{\theta n}=\sum_{\theta =0}^{2(p+1)-1} 
\left ( \sum_{k =-\infty }^{+\infty } x[k]\epsilon^{k\theta}\right ) \epsilon^{\theta n}.
\end{equation*}

Changing the order of the sums, the right side of
the expression above becomes
\begin{equation*}
\sum_{\theta =0}^{2(p+1)-1} X(\epsilon^\theta)\epsilon^{\theta n}=
\sum_{k =-\infty }^{+\infty } x[k] \left ( \sum_{\theta =0}^{2(p+1)-1} \epsilon^{\theta(n-k)}\right ).
\end{equation*}
But, from Lemma 1, the internal sum is nonzero only for $k = n$, so that
\begin{equation*}
\sum_{\theta =0}^{2(p+1)-1} X(\epsilon^\theta)\epsilon^{\theta n}=
2(p+1) x[k] 
\end{equation*}
and the result follows.\\

It is interesting to notice that, although the direct FFDTFT involves an infinite sum, being capable of handling infinite sequences, its inverse requires only a finite sum over the phase group $G_\theta$. the finite field discrete time Fourier transform introduced in this work satisfies most properties of the usual DTFT defined over the complex field $\mathbb{C}$, such as linearity, time shift, scaling and so on.\\

\textbf{Example 4.} The inverse FFDTFT of the plane spectrum $X(\epsilon^\theta) = 1$ is
\begin{equation*}
x[n]=\frac{1}{2(p+1)} \sum_{\theta=0}^{2(p+1)-1} \epsilon^\theta
\end{equation*}
and, from Lemma 1, follows $x[n] = \delta[n]$, as expected.

\section{Redefining the Finite Field Fourier Transform}
The Discrete Fourier Transform is a commonly used tool in Electrical Engineering. In a general setting, the DFT of a sequence $v = (v_i) \in \mathbb{E}$, is the sequence $V = (V_k) \in \mathbb{F}$ of elements
\begin{equation*}
V_k \triangleq \sum_{i=0}^{N-1} v_i W^{ik}
\end{equation*}
where $i,k=0,1,...,N-1$ and $W$ is an Nth root of unity in $\mathbb{F}$.\\

If $\mathbb{E}$ is the field $\mathbb{R}$ of real numbers and $\mathbb{F}=\mathbb{C}$  then, $W = (\exp-j2\pi/N)$ and we have the usual DFT. In this case the transformed vector is, in general, complex.\\

If $\mathbb{E}=GF(p)$ and $\mathbb{F}=GF(p^m)$, with $m \geq 1$, then $W = a$ is an element of multiplicative order $N$ of $GF(p^m)$. In this case, we have the Finite Field Fourier Transform. The FFFT definition with a kernel $W = a \in GF(p^m)$, makes the transformation to be a real one.\\

For the above, a definition of the FFFT analogous to the usual DFT, should use a complex kernel. In this case, we have not only a more appropriate version of the FFFT, but also a greater flexibility in the transform length.\\

\textit{\textbf{Definition 9}. Let $f = (f_0,f_1,...,f_{N-1})$ be a vector of length $N$ with components over $GF(q)$, where $q = p^r$. Then the vector $F = (F_0,F_1,...,F_{N-1})$, with components over $GL(q^m)$ given by}
\begin{equation*}
F_k \triangleq \sum_{i=0}^{N-1} f_i \zeta^{ik}
\end{equation*}
\textit{where $\zeta$ is an element of order $N$ in $GL(q^m)$, is the Finite Field Fourier Transform of $f$}.\\

This FFFT has the same properties as the one introduced by Pollard \cite{Pollard}. Indeed, the last one turns out to be a particular case of the Definition 12, when $\zeta = a + jb$ and $b = 0$.\\

\textit{\textbf{Proposition 6}. The Finite Field Fourier Transform in Definition 9 has lengths $N$, which divide $q^{2m}-1$}.

\textit{Proof}: The transform has length $N$, given by the order of the element $\zeta \in  GL(q^m)$. Since $|GL(q^m)| = q^{2m}-1$, the result follows.\\

Observe that, since $q^{2m}-1 = (q^m-1)(q^m+1)$, new lengths which are divisors of $q^m+1$ are now possible for the FFFT.

\section{Final Remarks}

This paper deals with discrete Fourier transforms defined over finite fields. Initially, some mathematical preliminaries were presented, which lead to the construction of complex numbers over a finite field. Afterwards, the problem of constructing finite field transforms capable to deal with infinite sequences defined over a Galois Field $GF(p)$, was approached. To deal with such sequences, the concept of “Cesa\`ro convergence” was used, and it was shown that periodic sequences over $GF(p)$ converge to the arithmetic mean of the Cesa\`ro sums of the sequence. As a direct consequence of this, a new transform, the discrete time Fourier transform over a finite field, was introduced and its inversion formula was presented.\\

The Fourier transform over a finite field was also considered and a new definition for it was proposed. This formulation, not only generalizes it, but better mimics DFT, since it contains a complex kernel. The new definition allows a greater flexibility in the choice of lengths to the transform.\\

The transforms here introduced are able to process finite and infinite sequences and, hence, are useful tools to work with FIR and IIR filters defined over finite algebraic structures.

\bibliographystyle{unsrt}

\begin{thebibliography}{1}

\bibitem{Oppenheim}
A. V. Oppenheim, R. W. Schafer e J. R. Buck, \emph{Discrete-Time Signal Processing}, Prentice-Hall, 1999.

\bibitem{Pollard}
J. M. Pollard, The Fast Fourier Transform in a Finite Field, \emph{Math. Comput.}, vol. 25, No. 114, pp. 365-374, Apr. 1971.

\bibitem{Campello}
R. M. Campello de Souza, H. M. de Oliveira and A. N. Kauffman, Trigonometry in Finite Fields and a New Hartley Transform, Proc. of the \emph{IEEE Int. Symp. on Info. Theory}, p.293, Cambridge, MA, Aug. 1998.

\bibitem{de_Oliveira1}
H. M. de Oliveira, R. M. Campello de Souza and A. N. Kauffman, Efficient Multiplex for Band-Limited Channels: Galois-Field Division Multiple Access, \emph{Proceedings of the 1999 Workshop on Coding and Cryptography} - WCC ’99, pp. 235-241, Paris, Jan. 1999.

\bibitem{de_Oliveira2}
H. M. de Oliveira, R. M. Campello de Souza, Orthogonal Multilevel Spreading Sequence Design, in \emph{Coding, Communications and Broadcasting}, pp. 291-303, Eds. P. Farrell, M. Darnell and B. Honary, Research Studies Press / John Wiley, 2000.

\bibitem{Cooklev}
T. Cooklev, A. Nishihara and M. Sablatash, Theory of Filter banks over Finite Fields, \emph{IEEE Asia Pacific Conference on Circuits and Systems}, APCCAS’94, pp. 260-265, 1994.

\bibitem{Blahut}
R. E. Blahut, \emph{Fast Algorithms for Digital Signal Processing}, Addison Wesley, 1985.

\bibitem{Lidl}
R. Lidl e H. Niederreiter, \emph{Introduction to Finite Fields and their Applications}, Cambridge University press, 1986.

\bibitem{Burton}
D.M. Burton, \emph{Elementary Number Theory}, McGraw-Hill, 1997.

\bibitem{Bartle}
R. G. Bartle, \emph{The Elements of Real Analysis}, John Wiley, 1967.

\end{thebibliography}

\end{document}